\documentclass{article}
\usepackage[utf8]{inputenc}
\usepackage{amsmath}
\usepackage[style=numeric]{biblatex}

\DeclareMathOperator{\diag}{diag}

\title{An even more straightforward proof of Descartes's circle theorem}
\author{Alden Bradford}
\date{September 2022}

\begin{document}

\maketitle
\begin{abstract}
How was this proof overlooked for 181 years? We give a simple proof of Descartes's circle theorem using Cayley-Menger determinants.
\end{abstract}

\section{Introduction}

Descartes's Circle Theorem states that the radii of four mutually tangent circles $r_1$, $r_2$, $r_3$, and $r_4$ satisfy
\[
\left(\frac{1}{r_1}+\frac{1}{r_2}+\frac{1}{r_3}+\frac{1}{r_4}\right)^2 = 2\left(\frac{1}{r_1^2}+\frac{1}{r_2^2}+\frac{1}{r_3^2}+\frac{1}{r_4^2}\right).
\]
The radii are chosen to be negative if the corresponding circle encloses the others. In this way we preserve the relation $d_{ij}^2=(r_i+r_j)^2$ where $d_{ij}$ is the distance between the centers of circles.

An article \cite{LevriePaul2019ASPo} published in this journal in 2019 gives a short history of the theorem and provides an original and straightforward proof based on Heron's formula. Here, we provide an even more straightforward proof based on a generalization of Heron's formula, the Cayley-Menger determinant.

\section{Cayley-Menger determinants}

The Cayley-Menger determinant was first introduced by Arthur Cayley in 1841 \cite{cdi_hathitrust_hathifiles_nyp_33433062744911, LibertiLeo2016Smgf}. It gives a formula for the volume of an $n$-simplex in terms of the pairwise distances between the vertices. In the case of a triangle with side lengths $a$, $b$, $c$ and area $A$ the formula is equivalent to Heron's formula,
\[
-16A^2 = 
\begin{vmatrix}
0 & 1 & 1 & 1\\
1 & 0 & c^2 & b^2\\
1 & c^2 & 0 & a^2\\
1 & b^2 & a^2 & 0
\end{vmatrix}.
\]
To prove Descartes's theorem, we will use the Cayley-Menger determinant for the tetrahedron,
\[
288v^2= \begin{vmatrix}
0&1&1&1&1\\
1&0&d_{12}^{2}&d_{13}^{2}&d_{14}^{2}\\
1&d_{21}^{2}&0&d_{23}^{2}&d_{24}^{2}\\
1&d_{31}^{2}&d_{32}^{2}&0&d_{34}^{2}\\
1&d_{42}^{2}&d_{42}^{2}&d_{43}^{2}&0\\
\end{vmatrix}.
\]
This formula is quite simple to show. We present the tetrahedral case here, adapted from the more general proof given in \cite{alma99130294800001081} for an $n$-simplex. Let $x_j$ be the $j$th vertex of the tetrahedron, and $x_{ij}$ its $i$th component. Write
\[
U = \begin{bmatrix}
1 & |x_1|^2 & |x_2|^2 & |x_3|^2 & |x_4|^2\\
0 & 1 & 1 & 1 & 1\\
0 & x_{11} & x_{12} & x_{13} & x_{14} \\
0 & x_{21} & x_{22} & x_{23} & x_{24} \\
0 & x_{31} & x_{32} & x_{33} & x_{34} 
\end{bmatrix},
\]
\[
W=\begin{bmatrix}
0 & 1 & 0 & 0 & 0\\
1 & 0 & 0 & 0 & 0\\
0 & 0 &-2 & 0 & 0\\
0 & 0 & 0 &-2 & 0\\
0 & 0 & 0 & 0 &-2
\end{bmatrix}.
\]
Because $d_{ij}^2 = |x_i-x_j|^2=|x_i|^2+|x_j|^2-2x_i\cdot x_j$ we have that $U^T W U = D$, the matrix of the Cayley-Menger determinant. On the other hand, $|W|=8$ and we can expand along the first column of $U$ to reach the standard cross-product style formula for the volume of a tetrahedron,
\[
|U|
= \begin{vmatrix}
1 & 1 & 1 & 1\\
x_{11} & x_{12} & x_{13} & x_{14} \\
x_{21} & x_{22} & x_{23} & x_{24} \\
x_{31} & x_{32} & x_{33} & x_{34} 
\end{vmatrix}
=
6v.
\]
Putting these together, $|D|=|U^TWU|=|U|^2|W|=288v^2$.

\section{Proof}

The strategy for this proof of Descartes's circle theorem is to consider the tetrahedron whose vertices are the centers of the given circles. We apply the Cayley-Menger determinant formula, replacing each $d_{ij}$  with $r_i+r_j$. After simplifying we reach the formula
\[
v^2 = \left(\frac{r_1r_2r_3r_4}{3}\right)^2\left[\left(\frac{1}{r_1}+\frac{1}{r_2}+\frac{1}{r_3}+\frac{1}{r_4}\right)^2 -2\left(\frac{1}{r_1^2}+\frac{1}{r_2^2}+\frac{1}{r_3^2}+\frac{1}{r_4^2}\right)\right].
\]
This equation gives us exactly what we need. If the four circle centers lie in a plane then their tetrahedron will have zero volume. Since none of the radii are zero, the term $r_1r_2r_3r_4$ must be nonzero and hence, 
\[
\left(\frac{1}{r_1}+\frac{1}{r_2}+\frac{1}{r_3}+\frac{1}{r_4}\right)^2 -2\left(\frac{1}{r_1^2}+\frac{1}{r_2^2}+\frac{1}{r_3^2}+\frac{1}{r_4^2}\right)=0.
\]
That is the whole of the proof. All that remains is to justify the volume formula given above.

Let $D$ be the matrix in the Cayley-Menger determinant formula. When we expand the terms $d_{ij}^2=r_i^2+2r_ir_j+r_j^2$ we are left with the term $r_i^2$ repeated along each row. We can eliminate it using the matrix 
\[
P = \begin{bmatrix}
1 & -r_1^2 & -r_2^2 & -r_3^2 & -r_4^2\\
0 & 1 & 0 & 0 & 0\\
0 & 0 & 1 & 0 & 0\\
0 & 0 & 0 & 1 & 0\\
0 & 0 & 0 & 0 & 1
\end{bmatrix},
\]
giving us
\[
P^T D P = \begin{bmatrix}
0 & 1 & 1 & 1 & 1\\
1 & -2r_1^2 & 2r_1r_2 & 2r_1r_3 & 2r_1r_4\\
1 & 2r_2r_1 & -2r_2^2 & 2r_2r_3 & 2r_2r_4\\
1 & 2r_3r_1 & 2r_3r_2 & -2r_3^2 & 2r_3r_4\\
1 & 2r_4r_1 & 2r_4r_2 & 2r_4r_3 & -2r_4^2
\end{bmatrix}.
\]
Each column and row has a common factor of $r_i$, so we can pull it out using the matrix $Q = \diag(1, 1/r_1, 1/r_2, 1/r_3, 1/r_4)$. Then
\[
Q^T P^T D P Q = \begin{bmatrix}
0&\frac{1}{r_1}&\frac{1}{r_2}&\frac{1}{r_3}&\frac{1}{r_4}\\
\frac{1}{r_1}&-2&2&2&2\\
\frac{1}{r_2}&2&-2&2&2\\
\frac{1}{r_3}&2&2&-2&2\\
\frac{1}{r_4}&2&2&2&-2
\end{bmatrix}.
\]
Write $R = \begin{bmatrix} 1/r_1 & 1/r_2 & 1/r_3 & 1/r_4 \end{bmatrix}^T$ and $S= 2\mathbf{1}\mathbf{1}^T - 4I$ where $\mathbf{1}=\begin{bmatrix} 1& 1 & 1 & 1\end{bmatrix}$, which allows us to put this in the block form
\[
Q^T P^T D P Q = \begin{bmatrix}
0 & R^T \\
R & S
\end{bmatrix}.
\]

There are a couple things to note about $S$. First, observe that $S^2 = 16I$ and so $S^{-1} = \frac{1}{16} S$. Second, we can readily compute $|S|=-256$. We take advantage of both of these when applying a common rule for the determinants of block matrices. Recall that if $A_{22}$ is invertible then
\[
\begin{vmatrix}
A_{11} & A_{12}\\
A_{21} & A_{22}
\end{vmatrix}
=|A_{22}||A_{11}-A_{12}A_{22}^{-1} A_{21}|.
\]
We apply this here to find
\begin{align*}
|Q^T P^T D P Q| 
&= -|S| R^T S^{-1} R\\
&= 16 R^T S R\\
&= 32 R^T[\mathbf{1}\mathbf{1}^T - 2I]R\\
&= 32 \left[ (R^T\mathbf{1})(\mathbf{1}^TR) - 2 R^T R\right]\\
&= 32 \left[ \left(\frac{1}{r_1}+\frac{1}{r_2}+\frac{1}{r_3}+\frac{1}{r_4}\right)^2 -2\left(\frac{1}{r_1^2}+\frac{1}{r_2^2}+\frac{1}{r_3^2}+\frac{1}{r_4^2}\right)\right].
\end{align*}
On the other hand, 
\begin{align*}
|Q^T P^T D P Q|
&= |P|^2 |Q|^2 |D|\\
&= (1)^2\left(\frac{1}{r_1r_2r_3r_4}\right)^2(288v^2).
\end{align*}
This completes the proof.

\section{Generalizing to higher dimensions}
The Soddy-Gosset theorem generalizes the Descartes circle theorem to configurations of $n+2$ pairwise-tangent spheres in $n$ dimensions. Taking $r_1$, $r_2$, ..., $r_{n+2}$ to be the signed  radii of the $n$-dimensional spheres, the theorem states that
\[
\left(\frac{1}{r_1}+\cdots+\frac{1}{r_{n+2}}\right)^2 = n\left(\frac{1}{r_1^2}+\cdots+\frac{1}{r_{n+2}^2}\right).
\]
The above proof generalizes perfectly well to higher dimensions, giving exactly the Soddy-Gosset theorem. The only changes are the sizes of the matrices. In general the Cayley-Menger determinant for $n+2$ points evaluates to 
\[(-1)^{n}2^{n+1}\left((n+1)!v_{n+1}\right)^2,
\]
where $v_{n+1}$ is written to emphasize that the formula gives an $n+1$-dimensional volume. The general matrix $S$ has determinant and inverse
\begin{align*}
|S|&=(-1)^{n+1}2^{2n+3}n, \\
S^{-1} &= \frac{1}{4n}\mathbf{1}\mathbf{1}^T-\frac{1}{4}I.
\end{align*}
Carrying these changes through the computation gives 
\[
v_{n+1}^2=2^n \left(\frac{r_1\cdots r_{n+2}}{(n+1)!}\right)^2 \left[\left(\frac{1}{r_1}+\cdots+\frac{1}{r_{n+2}}\right)^2-n\left(\frac{1}{r_1^2}+\cdots+\frac{1}{r_{n+2}^2}\right)\right].
\]

\section{Analysis}
The proof above has an intuitive meaning, as it connects the Descartes formula to the volume of a simplex. It generalizes nicely to higher dimensions as well. Since Descartes's theorem deals only with the distances between points it seems natural to approach this problem from the perspective of distance geometry. Cayley-Menger determinants are a foundational tool in distance geometry. In Cayley's original paper he used the argument that the determinant should be zero when all the points lie in a lower-dimensional subspace, which is the same way we use it here. The simplex which joins the centers in a sphere packing has been used before to prove facts about sphere packings \cite{CooperDaryl1996Csca, https://doi.org/10.48550/arxiv.math/0302069, ChenHao2016ABPa}. Certainly many brilliant minds have sat down to prove Descartes's theorem, bringing all manner of advanced techniques to bear.

Given all these hints littered throughout history, why wasn't this proof noticed 181 years ago when Cayley first introduced his determinant? We can only speculate of course, but I would suggest the reason is there was less interest in Descartes's circle theorem back then. We are in the middle of a renaissance of interest in the theorem, due both to the complex-valued generalization \cite{LagariasJeffreyC.2002BtDC} and interest in the number-theoretic properties of Apollonian circle packings \cite{SarnakPeter2011IAP}. The proof given here does not address the complex-valued generalization, which is the main limitation of this proof. There is a concise and elegant proof of that due to Kocik \cite{https://doi.org/10.48550/arxiv.1910.09174}.

\section{Acknowledgements}

I thank Edna Jones for answering many of my questions about Descartes's theorem and Apollonian circle packings. I also thank the reviewers for several helpful suggestions.

\printbibliography

@article{LevriePaul2019ASPo,
year = {2019},
author = {Levrie, Paul},
address = {New York},
copyright = {Springer Science+Business Media, LLC, part of Springer Nature 2019},
issn = {0343-6993},
journal = {The Mathematical intelligencer},
language = {eng},
number = {3},
pages = {24-27},
publisher = {Springer US},
title = {A Straightforward Proof of Descartes’s Circle Theorem},
volume = {41},
}

@article{LibertiLeo2016Smgf,
year = {2016},
author = {Liberti, Leo and Lavor, Carlile},
address = {Oxford},
copyright = {2015 The Authors. International Transactions in Operational Research © 2015 International Federation of Operational Research Societies Published by John Wiley & Sons Ltd, 9600 Garsington Road, Oxford OX4 2DQ, UK and 350 Main St, Malden, MA02148, USA.},
issn = {0969-6016},
journal = {International transactions in operational research},
language = {eng},
number = {5},
pages = {897-920},
publisher = {Blackwell Publishing Ltd},
title = {Six mathematical gems from the history of distance geometry},
volume = {23},
}

@misc{cdi_hathitrust_hathifiles_nyp_33433062744911,
year = {1841},
address = {England},
language = {eng},
publisher = {E. Johnson, 1841-1846},
journal = {The Cambridge mathematical journal. : v. 2 (Nov. 1839-May 1841)},
title = {On a Theorem in the Geometry of Position},
author = {Cayley, Arthur},
pages = {267-271}
}

@article{LagariasJeffreyC.2002BtDC,
year = {2002},
author = {Lagarias, Jeffrey C. and Mallows, Colin L. and Wilks, Allan R.},
issn = {0002-9890},
journal = {The American mathematical monthly},
language = {eng},
number = {4},
pages = {338-},
title = {Beyond the Descartes Circle Theorem},
volume = {109},
}

@article{SarnakPeter2011IAP,
year = {2011},
author = {Sarnak, Peter},
copyright = {Copyright the Mathematical Association of America 2011},
issn = {0002-9890},
journal = {The American mathematical monthly},
language = {eng},
number = {4},
pages = {291-306},
publisher = {Mathematical Association of America},
title = {Integral Apollonian Packings},
volume = {118},
}

@misc{https://doi.org/10.48550/arxiv.1910.09174,
  doi = {10.48550/ARXIV.1910.09174},
  
  url = {https://arxiv.org/abs/1910.09174},
  
  author = {Kocik, Jerzy},
  
  title = {Proof of Descartes circle formula and its generalization clarified},
  
  publisher = {arXiv},
  
  year = {2019},
  
  copyright = {arXiv.org perpetual, non-exclusive license}
}

@book{alma99130294800001081,
author = {Berger, Marcel},
address = {Berlin ;},
booktitle = {Geometry},
isbn = {0387116583},
language = {eng},
lccn = {87136730},
publisher = {Springer-Verlag},
series = {Universitext},
title = {Geometry },
year = {1987},
}

@article{CooperDaryl1996Csca,
author = {Cooper, Daryl and Rivin, Igor},
issn = {1073-2780},
journal = {Mathematical research letters},
language = {eng},
number = {1},
pages = {51-60},
title = {Combinatorial scalar curvature and rigidity of ball packings},
volume = {3},
year = {1996},
}

@misc{https://doi.org/10.48550/arxiv.math/0302069,
  doi = {10.48550/ARXIV.MATH/0302069},
  
  url = {https://arxiv.org/abs/math/0302069},
  
  author = {Rivin, Igor},
  
  title = {An extended correction to ``Combinatorial Scalar Curvature and Rigidity of Ball Packings,'' (by D. Cooper and I. Rivin)},
  
  publisher = {arXiv},
  
  year = {2003},
  
  copyright = {Assumed arXiv.org perpetual, non-exclusive license to distribute this article for submissions made before January 2004}
}

@article{ChenHao2016ABPa,
author = {Chen, Hao},
address = {New York},
copyright = {The Author(s) 2016},
issn = {0179-5376},
journal = {Discrete \& computational geometry},
language = {eng},
number = {4},
pages = {801-826},
publisher = {Springer US},
title = {Apollonian Ball Packings and Stacked Polytopes},
volume = {55},
year = {2016},
}

\end{document}